\documentclass[12pt]{article}
\usepackage{amscd}
\usepackage{amsfonts}
\usepackage{amsmath,amsthm,hyperref}
\usepackage[all]{xy}

\newtheorem{theorem}{Theorem}[section]
\newtheorem{proposition}[theorem]{Proposition}
\newtheorem{definition}[theorem]{Definition}

\theoremstyle{remark}
\newtheorem{remark}[theorem]{Remark}
\newtheorem{example}[theorem]{\bf Example}

\begin{document}
\title{\bf{Polar transform of Spacelike isothermic surfaces in 4-dimensional
Lorentzian space forms}}
\author{{Xiang Ma\footnote{
Project 10771005 supported by NSFC.} ~~and~ Peng Wang } \\[3mm]
{\small
Dedicated to Professor Udo Simon on the occasion of his 70th
birthday}}
\date{}
\maketitle

\begin{center}
{\bf Abstract}
\end{center}

The conformal geometry of spacelike surfaces in 4-dimensional
Lorentzian space forms has been studied by the authors in a
previous paper, where the so-called polar transform was
introduced. Here it is shown that this transform preserves
spacelike conformal isothermic surfaces. We relate this new
transform with the known transforms (Darboux transform and
spectral transform) of
isothermic surfaces by establishing the permutability theorems. \\


\section{Introduction}

Isothermic surfaces are classical objects in differential
geometry. The most beautiful results about them are those
transforms producing new isothermic surfaces, such as the dual
isothermic surface (also named the Christoffel transform), the
spectral transform (also known as the T-transform, the Bianchi
transform or the Calapso transform), and the Darboux transform. In
particular, people established the permutability theorems relating
them (see \cite{Jer} for an overview). These facts indicate that
there is a structure of integrable system underlying the theory
about isothermic surfaces, which was revealed only in the past 20
years \cite{CGS, BHPP, Bur, BDPT}.

For Lorentzian space forms there is also a parallel theory of
conformal geometry \cite{AP}. Thus it is natural to study
isothermic surfaces in this context \cite{Du, FI}. Zuo et al
\cite{ZCC} generalized the Darboux transform of isothermic
surfaces to the pseudo-Riemannian space forms using the methods
developed by Burstall in \cite{Bur} and Bruck et al in
\cite{BDPT}. Their methods mainly concerned the integrable system
aspect of the theory.

In \cite{Ma-W} we studied spacelike surfaces in $Q^4_1$, the
conformal compactification of the 4-dimensional Lorentzian space
forms $R^{4}_{1},S^{4}_{1}$ and $H^{4}_{1}$. The key observation
is that in this codim-2 case, the normal plane at any point is
Lorentzian. The two null lines $[L],[R]$ in this plane define two
conformal maps into $Q^4_1$, called the \emph{ left} and the
\emph{right polar surface}, respectively. Conversely, $Y$ is also
the right polar surface of $[L]$, and the left polar surface of
$[R]$ (when $[L]$ and $[R]$ are immersions). Applying these
transforms successively, we obtain a sequence of conformal
immersions. We proved in \cite{Ma-W} that these transforms
preserve the Willmore property.
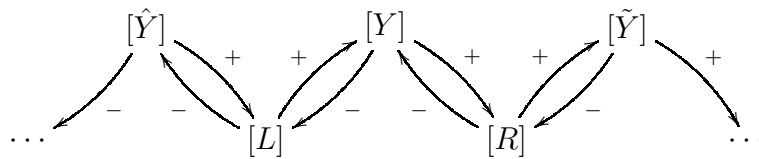
\[
\xymatrix{& &[\hat{Y}] \ar@/^/[dl]^{-} \ar@/^/[dr]^{+} & &[Y]
\ar@/^/[dl]^{-} \ar@/^/[dr]^{+}
& &[\tilde{Y}] \ar@/^/[dl]^{-} \ar@/^/[dr]^{+} &\\
&\cdots & &[L]\ar@/^/[ul]^{-} \ar@/^/[ur]^{+} & &[R]\ar@/^/[ul]^{-}
\ar@/^/[ur]^{+} & &\cdots}
\]

The first main result in this paper is that the isothermic
property is also invariant under the polar transform
(Theorem~\ref{thm-isoth}). A new isothermic surface produced in
this way is neither the spectral transform nor the Darboux
transform of $[Y]$. Hence it turns out to be a new transform for
isothermic surfaces. (The authors note that similar results hold
for timelike Willmore surfaces and for timelike isothermic
surfaces in $Q^4_2$, which might be treated in another paper.)

It is natural to wonder about the relationship between this new
transform and the old ones. In particular, does the polar
transform commute with the spectral transform or the Darboux
transform? The answer is affirmative. Two of such
\emph{permutability theorems} are established at here. See
Theorem~\ref{thm-commu1} and Theorem~\ref{thm-commu2}.

This paper is organized as follows. In Section~2 and Section~3 we
review the main theory about the Lorentzian conformal space
$Q^{4}_{1}$ and spacelike surfaces in it. The definition and
examples of isothermic surfaces are discussed in Section~4. Then
we introduce the polar transform of spacelike isothermic surfaces
in Section~5. Finally, after describing the spectral transform and
Darboux transform of an isothermic surface, we establish the
commutability between them and the polar transform in Section~6
and Section~7 separately.

\section{The Lorentzian conformal space $Q^{4}_{1}$}

Let $\mathbb{R}^n_s$ denote the space $\mathbb{R}^n$ equipped with
the quadratic form \[\langle
x,x\rangle=\sum^{n-s}_{1}x^2_i-\sum^n_{n-s+1}x^{2}_i.\] In this
paper we will mainly work with $\mathbb{R}^6_2$, whose light cone
is denoted as $C^5$. The quadric
\[
Q^4_1=\{\ [x]\in\mathbb{R}P^5\ |\ x\in C^5\setminus \{0\} \}
\]
is exactly the projective light cone with the projection map
$\pi:C^5\setminus\{0\}\rightarrow Q^4_1$. It is easy to see that
$Q^4_1$ is equipped with a Lorentzian metric $h$ induced from
projection $S^3\times S^1\rightarrow Q^4_1$. Here
\begin{equation}S^3\times
S^1=\{x\in\mathbb{R}^6_2\ |\
\sum^4_{i=1}x^{2}_{i}=x^2_5+x^{2}_6=1\}\subset C^5\setminus\{0\}
\end{equation}
is endowed with the Lorentzian metric $g(S^3)\oplus (-g(S^1))$,
where $g(S^3)$ and $g(S^1)$ are standard metrics on $S^3$ and
$S^1$. The conformal group of $(Q^4_1,[h])$ is exactly the
orthogonal group $O(4,2)/\{\pm1\}$, which keeps the inner product
of $\mathbb{R}^6_2$ invariant and acts on $Q^4_1$ by
\begin{equation}
T([x])=[xT],\ T\in O(4,2).
\end{equation}
As in Moebius geometry, $Q^4_1$ serves as the common conformal
compactification of the three 4-dimensional Lorentzian space forms
given below, each with constant sectional curvature $c=0,+1,-1$:
\begin{align*}
R^4_1 &, ~c=0; \\
S^4_1 &:=\{x\in\mathbb{R}^5_1\ |\ \langle x,x\rangle=1 \}, c=1; \\
H^4_1 &:=\{x\in\mathbb{R}^5_2\ |\ \langle x,x\rangle=-1 \}, c=-1.
\end{align*}
The conformal embedding into $Q^4_1$ for each of them is
\begin{equation}\label{varphi}
\begin{array}{llll}
\varphi_{0}:R^4_1\rightarrow Q^4_1, ~&
\varphi_{0}(x)=[(\frac{-1+\langle x,x\rangle}{2},
x,\frac{1+\langle x,x\rangle}{2})]; \\[1mm]
\varphi_{+}:S^4_1\rightarrow Q^4_1,
~& \varphi_{+}(x)=[(x,1)]; \\[1mm]
\varphi_{-}:H^4_1\rightarrow Q^4_1, ~& \varphi_{-}(x)=[(1,x)].
\end{array}
\end{equation}
Thus $Q^4_1$ is the proper space to study the conformal geometry
of these Lorentzian space forms.

We also have round spheres as the most important conformally
invariant objects in $Q^4_1$. Here we only discuss \emph{round
2-spheres} (they were named \emph{conformal 2-spheres} in
\cite{AP}). Each of them could be viewed as a geodesic 2-sphere in
a 3-dim Lorentzian space form. Alternatively, a round 2-sphere is
identified with a 4-dim Lorentzian subspace in $\mathbb{R}^6_2$.
Given such a 4-space $V$, the round 2-sphere is given by
\[S^2(V):=\{[v]\in Q^4_1 ~|~  v\in V\}.\]

\section{Spacelike surfaces in $Q^{4}_{1}$}

For a surface $y:M\rightarrow Q^{4}_{1}$ and any open subset
$U\subset M$, a local lift of $y$ is just a map $Y:U\rightarrow
C^5\setminus\{0\}$ such that $\pi\circ Y=y$. Two different local
lifts differ by a scaling, so the metric induced from them are
conformal to each other.

Let $M$ be a Riemann surface. An immersion $y:M\rightarrow
Q^{4}_{1}$ is called a conformal spacelike surface if $\langle
Y_{z},Y_{z}\rangle=0$ and $\langle Y_{z},Y_{\bar{z}}\rangle >0$
for any local lift $Y$ and any complex coordinate $z$ on $M$. For
such a surface there is a decomposition $M\times
\mathbb{R}^{6}_{2}=V\oplus V^{\perp}$, where
\begin{equation}
V={\rm Span}\{Y,{\rm Re}(Y_z),{\rm Im}(Y_z),Y_{z\bar{z}}\}
\end{equation}
is a Lorentzian rank-4 subbundle independent to the choice of $Y$
and $z$. $V^{\perp}$ is also a Lorentzian subbundle, which might be
identified with the normal bundle of $y$ in $Q^{4}_{1}$. Their
complexifications are denoted separately as $V_{\mathbb{C}}$ and
$V^{\perp}_{\mathbb{C}}$.

Fix a local coordinate $z$. There is a local lift $Y$ satisfying
$|{\rm d}Y|^2=|{\rm d}z|^2$, called the canonical lift (with respect
to $z$). Choose a frame $\{Y,Y_{z},Y_{\bar{z}},N\}$ of
$V_{\mathbb{C}}$, where $N\in\Gamma(V)$ is uniquely determined by
\begin{equation}\label{eq-N}
\langle N,Y_{z}\rangle=\langle N,Y_{\bar{z}}\rangle=\langle
N,N\rangle=0,\langle N,Y\rangle=-1.
\end{equation}
For $V^{\perp}$ which is a Lorentzian plane at every point of $M$, a
natural frame is $\{L,R\}$ such that
\begin{equation}\label{eq-LR}
\langle L,L\rangle=\langle R,R\rangle=0,\langle L,R\rangle=-1.
\end{equation}
Given frames as above, we note that $Y_{zz}$ is orthogonal to $Y$,
$Y_{z}$ and $Y_{\bar{z}}$. So there must be a complex function $s$
and a section $\kappa\in \Gamma(V_{\mathbb{C}}^{\perp})$ such that
\begin{equation}
Y_{zz}=-\frac{s}{2}Y+\kappa.
\end{equation}
This defines two basic invariants $\kappa$ and $s$ dependent on
$z$. Similar to the case in M\"obius geometry, $\kappa$ and $s$
are called the \emph{conformal Hopf differential} and the
\emph{Schwarzian derivative} of $y$, respectively (see
\cite{BPP},\cite{Ma1}). Decompose $\kappa$ as
\begin{equation}\label{eq-kappa}
\kappa=\lambda_{1}L+\lambda_{2}R.
\end{equation}
Let $D$ denote the normal connection, i.e. the induced connection
on the bundle $V^{\perp}$. We have
\[D_z L=\alpha L,~~D_z R=-\alpha R \]
for the connection 1-form $\alpha{\rm d}z$. Denote
\begin{equation}\label{eq-Dkappa}
\langle \kappa,\bar\kappa\rangle=-\beta,~~ D_{\bar z}\kappa=\gamma_1
L +\gamma_2 R,
\end{equation}
where
\begin{equation}\label{eq-gamma}
\left\{\begin{array}{llll}
\beta=\lambda_1\bar{\lambda}_2+\lambda_2\bar{\lambda}_1,\\
\gamma_1=\lambda_{1\bar{z}}+\lambda_1\bar{\alpha},\\
\gamma_2=\lambda_{2\bar{z}}-\lambda_2\bar{\alpha}.
\end{array}\right.
\end{equation}
The structure equations are given as follows:
\begin{equation}\label{eq-moving}
\left\{\begin {array}{lllll}
Y_{zz}=-\frac{s}{2}Y+\lambda_{1}L+\lambda_{2}R,\\
Y_{z\bar{z}}=\beta Y+\frac{1}{2}N,\\
N_{z}=2\beta Y_{z}-sY_{\bar{z}}+2\gamma_{1}L+2\gamma_{2}R,\\
L_{z}=\alpha L-2\gamma_{2}Y+2\lambda_{2}Y_{\bar{z}},\\
R_{z}=-\alpha R-2\gamma_{1}Y+2\lambda_{1}Y_{\bar{z}},
\end {array}\right.
\end{equation}
The conformal Gauss, Codazzi and Ricci equations as integrable
conditions are:
\begin{equation}\label{eq-integ}
\left\{\begin {array}{lllll}
s_{\bar{z}}=-2\beta_{z}-4\lambda_{1}\bar{\gamma}_{2}
-4\lambda_{2}\bar{\gamma}_{1},\\
{\rm Im}(\gamma_{1\bar{z}}+\gamma_{1}\bar{\alpha}
+\frac{\bar{s}}{2}\lambda_{1})=0,\\
{\rm Im}(\gamma_{2\bar{z}}-\gamma_{2}\bar{\alpha}
+\frac{\bar{s}}{2}\lambda_{2})=0,\\
\alpha_{\bar{z}}-\bar\alpha_z
=2(\lambda_{1}\bar{\lambda}_{2}-\lambda_{2}\bar{\lambda}_{1}).
\end {array}\right.
\end{equation}

\section{Spacelike isothermic surfaces}

\begin{definition}
Let $y:M\rightarrow Q^{4}_{1}$ be a conformal spacelike surface
without umbilic points. It is called isothermic if around each
point of $M$ there exists a complex coordinate $z$ and canonical
lift $Y$ such that the Hopf differential $\kappa$ is real-valued.
Such a coordinate $z$ is called an adapted coordinate.
\end{definition}
Since $\kappa$ is real-valued, from the conformal Ricci equations
in \eqref{eq-integ} we see that its normal bundle is flat. This is
an important property of isothermic surfaces, which guarantees
that all shape operators commute and the curvature lines could
still be defined. Indeed we can equivalently define $y$ to be
isothermic if it has flat normal bundle and if it has conformal
curvature line parameters. Put differently, the two fundamental
forms of an isothermic surface are of the form
\begin{equation}\label{iso}
I=e^{2\omega}(du^2+dv^2),\ II=(b_1 du^2+b_2 dv^2)e_{3}+(b_3
du^2+b_4 dv^2)e_{4}
\end{equation}
with respect to some parallel normal frame $\{e_3,e_4\}$. Then
$(u,v)$ are curvature line parameters and $z=u+{\rm i}v$ is an
adapted complex coordinate.

Our definition generalizes the notion of isothermic surfaces in
3-dim space forms and includes them as special cases. In the
following we provide more examples of isothermic surfaces in
$Q^4_1$.
\begin{example}
Rotational surfaces in $\mathbb{R}^3$ are isothermic as well
known. To generalize this construction, consider a spacelike curve
$\gamma(u)=(0,f(u),g(u),h(u)):\mathbb{R}\rightarrow
\mathbb{R}^{4}_{1}$ such that $f(u)\neq0,f'(u)\neq0,g'(u)\neq
h'(u)$. A rotational surface
$x:\mathbb{R}\times[0,2\pi]\rightarrow R^{4}_{1}$ generated by
$\gamma$ is just
\[x(u,v)=\Big(f(u)\cos
v,f(u)\sin v,g(u),h(u)\Big).\] It is easy to verify that
\eqref{iso} is satisfied when $u$ is reparameterized suitably.
\end{example}
\begin{example}
In \cite{Ma-W} we constructed a class of homogenous spacelike tori
as below, which are both Willmore and isothermic. Set
$\psi=\psi(t,\theta)=\theta/\sqrt{t^{2}-1}$. Then
$Y_{t}(\theta,\phi):\mathbb{R}\times\mathbb{R}\rightarrow
\mathbb{R}^{6}_{2}$ is given by
\[
Y_{t}(\theta,\phi)=\Big(\cos(t\psi)\cos\phi,\cos(t\psi)\sin\phi,
\sin(t\psi)\cos\phi,\sin(t\psi)\sin\phi,\cos\psi,\sin\psi\Big).
\]
Note that the period condition is satisfied if $t$ is a rational
number; hence after projection $\pi$ we obtain an immersed torus.
For the details see \cite{Ma-W}.
\end{example}

\section{Polar transform of isothermic surfaces}

For a conformal spacelike surface $y:M\rightarrow Q^4_1$ with
canonical lift $Y:M\rightarrow \mathbb{R}^6_2$ with respect to
complex coordinate $z=u+{\rm i}v$, its normal plane at any point
is spanned by two lightlike vectors $L,R$, determined up to a real
factor around each point. Suppose that $\mathbb{R}^6_2$ is endowed
with a fixed orientation and that
\[
\{Y,Y_u,Y_v,N,R,L\}
\]
form a positively oriented frame. $\{R,L\}$ might also be viewed as
a frame of the normal plane compatible with the orientation of $M $
and that of the ambient space. Since $\langle L,R\rangle=-1$ has
been fixed in \eqref{eq-LR}, either one of the null lines $[L]$
($[R]$) is well-defined.
\begin{definition}
The two maps $[L],[R]:M^2 \rightarrow Q^{4}_{1}$ are named the
\emph{left} and the \emph{right polar surface} of $y=[Y]$,
respectively.
\end{definition}

Alternatively sometimes we call $[L],[R]$ the (left and right)
polar transforms of $[Y]$. An interesting fact showed in
\cite{Ma-W} is that they again produce conformal mappings;
Moreover we have a \emph{duality} in this construction:

\begin{proposition}[\cite{Ma-W}]\label{prop-inverse}
The polar surfaces $[L],[R]:M^2\to Q^{4}_{1}$ are both conformal
maps. $[L]$ $([R]$) is degenerate if, and only if, $\lambda_2=0$
$(\lambda_1=0$); it is a spacelike immersion otherwise. The
original surface $[Y]$ is the left polar surface of $[R]$ $($the
right polar surface of $[L])$ when $[R]$ $([L])$ is not
degenerate.
\end{proposition}
In \cite{Ma-W} we have shown that the polar transforms of a
spacelike Willmore surface are again Willmore. Here we want to
show that a similar result holds true for isothermic surfaces.
\begin{theorem}\label{thm-isoth}
Let $y:M\rightarrow Q^{4}_{1}$ be a spacelike isothermic surface.
Then its left and right polar surfaces $[L],[R]:M\rightarrow
Q^{4}_{1}$ are also spacelike isothermic surfaces when they are
not degenerate. In particular they share the same adapted
coordinate $z$.
\end{theorem}
\begin{proof}
Let $Y:M \rightarrow \mathbb{R}^6_2$ be the canonical lift and
$\kappa$ be the real-valued conformal Hopf differential for an
adapted isothermic coordinate $z$. We show the conclusion for
$[L]$. For $[R]$ the proof is similar.

Assume that the left polar surface $[L]$ is an immersion, i.e.
$\lambda_{2}\ne 0$. Choosing $L$ such that $\kappa=\lambda_1
L+\lambda_2 R$ with $\lambda_{2}=\frac{1}{2}$, by
\eqref{eq-moving} we have
\begin{equation}\label{eq-lz}
L_{z}=\alpha L+\bar{\alpha}Y+Y_{\bar{z}}.
\end{equation}
Thus $L$ is the canonical life of $[L]:M\to Q^4_1$ as desired. To
determine the normal bundle of $[L]$, we differentiate once more
and invoke \eqref{eq-moving}, obtaining
\begin{equation}
L_{z\bar{z}}=(\alpha_{\bar{z}}+\lambda_1)L+\frac{1}{2}\left[R+2\alpha
Y_{z}+2\bar{\alpha}Y_{\bar{z}}+2|\alpha|^2L+2(\bar{\alpha}_{\bar{z}}
+\alpha^{2}-\frac{\bar{s}}{2})Y\right].
\end{equation}
We point out that each of $\alpha_{\bar{z}},\lambda_1$ and
$\alpha_z -\alpha^2-\frac{s}{2}$ is real valued (or by the Codazzi
and Ricci equations \eqref{eq-integ}). Now we can verify directly
that $Y$ and
\begin{equation}\label{eq-yhat}
\hat{Y}=N+2\alpha Y_{\bar{z}}+2\bar{\alpha}Y_{z}+2|\alpha|^2 Y
-2(\alpha_{z}-\alpha^{2}-\frac{s}{2})L
\end{equation}
are two lightlike vectors in the orthogonal complement of ${\rm
Span}\{L,L_u,L_v,L_{z\bar{z}}\}$ with $\langle
Y,\hat{Y}\rangle=-1$. Differentiate \eqref{eq-lz} at both sides.
After simplification we get
\begin{equation}\label{eq-lzz}
L_{zz}=(2\alpha_{z}-\frac{s}{2})L+\frac{1}{2}\hat{Y}
+(\bar{\alpha}_{z}+\lambda_1)Y.
\end{equation}
By definition, the conformal Hopf differential of $L$ is given by
$\kappa_{L}=-\frac{1}{2}\hat{Y}-(\bar{\alpha}_{z}+\lambda_1)Y$,
which is obviously real-valued. This shows that $[L]$ is
isothermic with the same adapted coordinate.
\end{proof}
Note that $\{L,L_u,L_v,L_{z\bar{z}},Y,\hat{Y}\}$ is again a
positively oriented frame. So $[Y]$ and $[\hat{Y}]$ is the right
and the left polar surface of $[L]$, respectively. This proves the
conclusion of Proposition~\ref{prop-inverse} in this special case.
On the other hand, $[\hat{Y}]$ is the left polar surface of $[L]$,
hence the 2-step left polar transform of $[Y]$.

\section{Permutability with spectral transform}

Let $y:M \rightarrow Q^{4}_{1}$ be an immersed spacelike
isothermic surface with canonical lift $Y:M \rightarrow
\mathbb{R}^6_2$ with respect to an adapted coordinate $z$. The
conformal Gauss, Codazzi, and Ricci equations are still satisfied
under the deformation
\[
s^{c}=s+c,~\lambda_1^c=\lambda_1,~\lambda_2^c=\lambda_2,
~\alpha^c=\alpha,
\]
where $c\in\mathbb{R}$ is a real parameter. By
the integrable conditions, there are an associated family of
non-congruent isothermic surfaces $[Y^c]$ with corresponding
invariants. Similar to the case of M\"obius geometry, they are
called the \emph{spectral transforms} of the original surface (see
\cite{BPP}). Observe that they are conformal and share the same
adapted coordinate $z$.

Now we have two transforms, the polar transform and the spectral
transform, associated with an isothermic surface. The
permutability between them is established as below.
\begin{theorem}\label{thm-commu1}
Let $y^c$ be a spectral transform (with parameter $c$) of
$y:M\rightarrow Q^{4}_{1}$, both being spacelike isothermic
surfaces. Denote their canonical lift as $[Y], [Y^c]$ for the same
adapted coordinate $z$. If the left polar surface $[L]$ and
$[L^c]$ corresponding to them are non-degenerate, then $[L^c]$ is
also a spectral transform (with parameter $c$) of $[L]$, i.e., we
have the commuting diagram:
\begin{equation*}
\begin{CD}
[Y]@>{}>> [Y^c] \\
@V{}VV @V{}VV\\
[L] @>{}>> ~[L^c]\
\end{CD}
\end{equation*}
A similar result holds between the right polar transform and the
spectral transform.
\end{theorem}
\begin{proof}
Set $Y_{zz}=-\frac{s}{2}Y+\lambda_{1}L+\lambda_{2}R$, and
$D_{z}L=\alpha L$. Choose $L$ such that $\lambda_{2}=\frac{1}{2}$.
By assumption, for $[Y^c]$ the corresponding frame
$\{Y^c,Y^c_{z},Y^c_{\bar{z}},N^c,L^c,R^c\}$ has the same inner
product matrix and satisfies
\[
Y^c_{zz}=-\frac{s+c}{2}Y^c+\lambda_{1}L^c+\frac{1}{2}R^c~,\ \
D^c_zL^c=\alpha L^c.
\]
Recall that we have computed out \eqref{eq-lzz}\eqref{eq-yhat}:
\begin{align*}
L_{zz}&=-(\frac{s}{2}-2\alpha_{z}
)L+\frac{1}{2}\hat{Y}+(\bar{\alpha}_{z}+\lambda_1)Y ,\\
\hat{Y}&=N+2\alpha Y_{\bar{z}}+2\bar{\alpha}Y_{z}+2|\alpha|^2Y
-2(\alpha_{z}-\alpha^{2}-\frac{s}{2})L,
\end{align*}
where $\{Y,\hat{Y}\}$ form a basis of the normal plane of $L$ at
any point. The same result applys to $Y^c$ and $L^c$, hence
\[
L^c_{zz}=-(\frac{s+c}{2}-2\alpha_{z}
)L^c+\frac{1}{2}\hat{Y}^c+(\bar{\alpha}_{z}+\lambda_1)Y^c ,
\]
where $\hat{Y}^c,Y^c$ span the normal bundle of $[L^c]$.
Comparison shows that the Schwarzian derivative of $[L]$ is
$s-4\alpha_{z}$, while that of $[L^c]$ differs from it by $c$ as
we expected. Their conformal Hopf differential has the same
components $\frac12$ and $\bar{\alpha}_{z}+\lambda_1$. Finally,
the normal connection of $[L]$ is given by $\langle
Y_z,\hat{Y}\rangle =\alpha$, exactly the same as $[Y]$. So $[L^c]$
also share the same normal connection as $[Y^c]$, which is again
$\alpha$. We conclude that $[L^c]$ is exactly a spectral transform
of $[L]$ with parameter $c$.
\end{proof}
\begin{remark}
For a spacelike Willmore surface there is also an associated
family of Willmore surfaces, called the Willmore spectral
transform. One could show that this transform commutes with the
left/right polar transform. We did not notice this result in
\cite{Ma-W}, yet the proof is similar and easy.
\end{remark}

\section{Permutability with Darboux transform}

The most important transform of isothermic surfaces in
$\mathbb{R}^n$ is the Darboux transform. It is a second isothermic
surface obtained from the original one by integration, depending
on the choice of initial values and a real parameter. (So there
are many of them.) For such a pair of isothermic surfaces forming
Darboux transform to each other, a geometric characterization is
that they envelop one and the same 2-sphere congruence at
corresponding points, and that their conformal curvature lines are
preserved by this correspondence (see \cite{Bur}, \cite{Jer},
\cite{Ma1}). This description is easy to adapted to our case:
\\
\\
\noindent {\bf Definition and Proposition} {\it Let
$y:M\rightarrow Q^{4}_{1}$ denote a spacelike isothermic surface
with canonical lift $[Y]$ with respect to the adapted coordinate
$z$. A spacelike immersion $y^{\ast}: M\rightarrow Q^{4}_{1}$ is
called a {\rm Darboux transform} of $y$ if its local lift
$Y^{\ast}$ satisfies
\begin{equation}\label{eq-isoth}
Y^{\ast}_{z}\in {\rm Span}_{\mathbb{C}}\{Y^{\ast},Y,Y_{\bar{z}}\}.
\end{equation}
Note that this is well-defined, where $Y^{\ast}$ is not
necessarily the canonical lift. We have the following conclusions:
\par 1) $y,y^{\ast}$ are conformal; they envelop one and the same
round 2-sphere congruence given by ${\rm Span} \{Y,Y^{\ast},{\rm
d}Y)\}={\rm Span} \{Y,Y^{\ast},{\rm d}Y^{\ast}\}$.
\par 2) Set $\langle Y,Y^{\ast}\rangle =-1$. Then
$Y^{\ast}_z=\frac{\mu}{2}Y^{\ast}+\theta(Y_{\bar{z}}+\frac{\bar\mu}{2}Y)$,
where $\theta$ is a non-zero real constant. This Darboux transform
is specified as \emph{$D^{\theta}$-transform}.
\par 3) $Y^{\ast}$ is an isothermic surface sharing the same
adapted coordinate $z$. Hence the curvature lines of
~$y,\tilde{y}$ do correspond.}
\begin{proof}
The conclusion 1) is obvious under the assumption
\eqref{eq-isoth}. (Recall that a round 2-sphere in
$\mathbb{R}^{4}_{1}$ is identified to a 4-dimensional Lorentzian
subspace in $\mathbb{R}^6_2$. See Section~2.)

The normalization $\langle Y,Y^{\ast}\rangle=-1$ ensures $\langle
Y_z,Y^{\ast}\rangle=-\langle Y^{\ast}_z,Y\rangle=\mu/2$. Then
$Y^{\ast}_{z}\in {\rm
Span}_{\mathbb{C}}\{Y^{\ast},Y,Y_{\bar{z}}\}$ is explicitly
expressed by
\begin{equation}\label{eq-touch}
Y^{\ast}_z=\frac{\mu}{2}Y^{\ast}+
\theta\left(Y_{\bar{z}}+\frac{\bar\mu}{2}Y\right ).
\end{equation}
Differentiate \eqref{eq-touch}. We obtain
\[
Y^{\ast}_{z\bar{z}}=\left(\frac{\mu}{2}\bar{\theta}\right)Y_z+
\left(\frac{\bar{\mu}}{2}\theta+\theta_{\bar{z}}\right)Y_{\bar{z}}
+\theta\bar{\kappa}
+\left(\frac{\mu_{\bar{z}}}{2}+\frac{\mu\bar{\mu}}{4}\right)Y^{\ast}
+(\cdots)Y.
\]
Since $\kappa$ is real and non-zero by assumption, comparison
shows that $\theta$ is real-valued with $\theta_{\bar{z}}=0$.
Hence $\theta$ must be a real constant. It is non-zero since
$[Y^{\ast}]$ is an immersion. This verifies 2). Note that
$\mu_{\bar{z}}$ is real by the same argument.

Now $\frac{1}{\theta}Y^{\ast}$ is the canonical lift of
$[Y^{\ast}]$. To show conclusion 3) we need only to show that
$Y^{\ast}_{zz}$ is real modulo the components of
$Y^{\ast},Y^{\ast}_z,Y^{\ast}_{\bar{z}},Y^{\ast}_{z\bar{z}}$.
Differentiate \eqref{eq-touch}. The result is
\begin{align*}
Y^{\ast}_{zz}&=\theta(Y_{z\bar{z}}
+\frac{\bar\mu}{2}Y_z+\frac{\bar\mu}{2}Y)
+(\cdots)Y^{\ast}+(\cdots)Y^{\ast}_z\\
&=\theta Y_{z\bar{z}}+\left (\frac{\mu_{\bar{z}}}{2}
-\frac{\mu\bar\mu}{4}\right )Y~~~~({\rm
mod}~~Y^{\ast},Y^{\ast}_z,Y^{\ast}_{\bar{z}}),
\end{align*}
which is real as desired. Finally, for $z=u+{\rm i}v$ with $u,v$
real, the $u$-curves and $v$-curves are exactly the curvature
lines on both of $y$ and $y^{\ast}$. This completes the proof.
\end{proof}

Write out $Y^{\ast}$ explicitly:
\[
Y^{\ast}=N+\bar{\mu}Y_z +\mu
Y_{\bar{z}}+(\frac{1}{2}|\mu|^2-4f_{1}f_{2})Y+2f_{1}L+2f_{2}R.
\]
Set
\begin{equation} P=Y_z+\frac{\mu}{2}Y,~\xi=L-2f_2Y,~\eta=R-2f_1Y.
\end{equation}
Note that both of $\xi,\eta$ are lightlike and $\langle
\xi,\eta\rangle=-1$. The orthogonal complement of ${\rm
Span}\{\xi,\eta\}$ gives the Ribaucour 2-sphere congruence
enveloped by $[Y],[Y^{\ast}]$. The structure equations
\eqref{eq-moving} of $Y$ can be rewritten with respect to the
frame $\{Y,Y^{\ast},P,\bar{P},\xi,\eta\}$ as below.
\begin{equation}\label{eq-moving-space-iso-1}
\left\{\begin {array}{lllll}
Y_{z}=-\frac{\mu}{2}Y+P,\\[1mm]
Y^{\ast}_{z}=~\frac{\mu}{2}Y^{\ast}+\theta \bar{P},\\[1mm]
P_z=~\frac{\mu}{2}P+\frac{\theta}{2}Y+\lambda_1\xi+\lambda_2\eta,\\[1mm]
\bar{P}_{z}=-\frac{\mu}{2}\bar{P}+\frac{1}{2}Y^{\ast}-f_1\xi-f_2\eta,\\[1mm]
\xi_{z}=-2f_{2}P+2\lambda_{2}\bar{P},\\[1mm]
\eta_{z}=-2f_{1}P+2\lambda_{1}\bar{P}.
\end {array}\right.
\end{equation}
Now let us find out the left polar transform of $Y^{\ast}$. From
\eqref{eq-moving-space-iso-1}, we have
\[
Y^{\ast}_{z\bar{z}}=\frac{\mu_{\bar{z}}}{2}Y^{\ast}
+\frac{\mu}{2}Y^{\ast}_{\bar{z}}+\theta \bar{P}_{\bar{z}}
=2f_1f_2Y^{\ast}+\frac{\theta^{2}}{2}N^{\ast},
\]
where
\[
N^{\ast}=Y+\frac{\mu}{\theta}P+\frac{\bar{\mu}}{\theta}\bar{P}
+\frac{|\mu|^2}{2\theta^2}Y^{\ast}
+\frac{2\lambda_1}{\theta}\xi+\frac{2\lambda_2}{\theta}\eta
-\frac{4\lambda_1\lambda_2}{\theta}Y^{\ast}.
\]
Set $L^{\ast}=\xi-\frac{2\lambda_2}{\theta}Y^{\ast},
R^{\ast}=\eta-\frac{2\lambda_1}{\theta}Y^{\ast}$. Plus the
orientation restriction, it is easy to see that
$[L^{\ast}],[R^{\ast}]$ are just the left and right polar
transform of $Y^{\ast}$. Suppose that $[L],[L^{\ast}]$ are both
non-degenerate. Computation shows
\begin{equation}
\begin{split}
L^{\ast}_{z}=-\frac{f_2}{\lambda_2}L_{\bar{z}}
-\frac{2\lambda_{2z}+\mu\lambda_2}{\lambda_2}(L-L^{\ast}).
 \end{split}
 \end{equation}
Note that by Theorem~\ref{thm-isoth} and \eqref{eq-moving}, $[L]$
has the same adapted coordinate $z$ with canonical lift
$\dfrac{1}{2\lambda_2}L$. Next, $-\dfrac{\theta}{2f_2}L^{\ast}$ is
a lift of $[L^{\ast}]$ such that
\begin{gather*}
\left\langle~\frac{1}{2\lambda_2}L,~-\frac{\theta}{2f_2}L^{\ast}~\right\rangle
=-\frac{1}{2\lambda_2}\frac{\theta}{2f_2} \left\langle ~\xi+2f_2
Y,~
\xi-\frac{2\lambda_2}{\theta}Y^{\ast}~\right\rangle=-1,\\
\left(-\frac{\theta}{2f_2}L^{\ast}\right)_z
=\theta\cdot\left(\frac{1}{2\lambda_2}L\right)_{\bar{z}}
+(\cdots)L+(\cdots)L^{\ast}.~~~~
\end{gather*}
This proves that $[L^{\ast}]$ is just a Darboux transform of $[L]$
with the same parameter $\theta$. So we have estalished
\begin{theorem}\label{thm-commu2}
Let $y:M\rightarrow Q^{4}_{1}$ be a spacelike isothermic surface
and $[y^{\ast}]$ be a $D^\theta$-transform of $y$. If both of
their left polar surfaces $[L]$ and $[L^{\ast}]$ are not
degenerate, $[L^{\ast}]$ is also a $D^\theta$-transform of $L$,
i.e., we have the commuting diagram:\begin{equation}
\begin{CD}
[Y] @>D^\theta-{\rm transform}>> [Y^{\ast}] \\
@V{left\ polar}VV @VV{left\ polar}V\\
[L] @>D^\theta-{\rm transform}>> ~[L^{\ast}]\
\end{CD}
\end{equation}
A similar result holds between the right polar transform and the
$D^\theta$-transform.
\end{theorem}

\def\refname{Reference}

\vspace{10mm} \noindent Xiang Ma, {\small\it LMAM, School of
Mathematical Sciences, Peking University, 100871 Beijing, People's
Republic of China. e-mail: {\sf maxiang@math.pku.edu.cn}}

\vspace{5mm}  \noindent Peng Wang, {\small\it LMAM, School of
Mathematical Sciences, Peking University, 100871 Beijing, People's
Republic of China. e-mail: {\sf wangpeng@math.pku.edu.cn}}

\end{document}